\documentclass[12pt]{article}
\usepackage{amsmath}
\usepackage{amsfonts}
\usepackage{amssymb}
\usepackage{amsthm}

\advance\textwidth20mm
\advance\textheight30mm
\advance\voffset-20mm
\advance\hoffset-15mm

\newtheorem{theor}{Theorem} %[section]
\newtheorem{lem}{Lemma} %[section]
\newtheorem{corol}{Corollary} %[section]

\newcommand{\dok}{{\it Proof.} \hspace{2mm}}
\newcommand{\edok}{$\square$\vspace{2mm}}

\begin{document}

\title{On complexity of %Combinatorial Complexity of %for the Number of edges in
envelopes of piecewise linear functions, unions and intersections
of polygons
}

\author{Pavel Kozhevnikov\footnotemark[1]
%\\ Moscow Institute of Physics and Technology.
}

\begingroup
\renewcommand{\thefootnote}{\fnsymbol{footnote}}

%\footnotetext[1]{e-mail: ilya.i.bogdanov@gmail.com \\
%Moscow Institute of Physics and Technology}
%\footnotetext[2]{e-mail: ivlevfyodor@mail.ru \\
% Moscow State University}
\footnotetext[1]{e-mail: p.kozhevn@gmail.com
\\ Moscow Institute of Physics and Technology}
\endgroup

%\footnote{Partially supported by
%a grant "AVCP Development of scientific potential of higher school", project 2.1.1/12136}}

\date{}

\maketitle

\begin{center}
{\bf Abstract}
\end{center}

%\vspace{1mm}

We prove tight upper bounds for the number of vertices of a simple
polygon that is the union or the intersection of two simple polygons with given numbers of
convex and concave vertices. The similar question on graphs of the lower (or upper) envelope
of two continuous piecewise linear functions is considered.

\vspace{3mm}


\vspace{3mm}


\subsection{Introduction %and presenting of results
}

\smallskip

%Опр. vertices
%Picewise-linear function define vertices $a_1< a_2< \ldots< a_n$
%We call the vertex {\it convex} if...
%%{\it concave}
%{\it reflex} if ....

%For a continuous piecewise-linear function $f_i:\mathbb{R}\to \mathbb{R}$
%consider its {\it vertices}. By a {\it convex} vertex (or
%a {\it reflex} vertex) of $f$ here we mean a convex (or reflex) vertex
%of  the set $\{(x, y)\, | \, y\geq f(x)\}$.

For a collection of functions $f_i : \mathbb{R} \to \mathbb{R}$, $i=1, 2, \ldots, k$,
the pointwise minimum function $f_0$  defined as
$f_0(x ) = \min \{f_1(x), \ldots, f_k(x)\}$ is called {\it the lower envelope} of %for  collection
$f_1, \ldots, f_k$. Similarly, {\it the upper envelope} defined as the pointwise maximum.

The epigraph of the lower envelope of $f_1, \ldots, f_k$ is the union of
the epigraphs of $f_1, \ldots, f_k$, while its hypograph is the intersection of
the hypographs of $f_1, \ldots, f_k$. That is the reason for considering envelopes
of piecewise-linear functions along with unions and intersection of polygons.

In this paper we deal with arrangements of two polygons or two graphs of
piecewise-linear functions
and focus on combinatorial complexity %(that is the number of vertices)
of a special face that is the common exterior (interior) of polygons or epigraph (hypograph)
of the envelopes.
Arrangements of polygons and graphs of
piecewise-linear functions (of one-variable) can be considered as
a special case of arrangement of segments and rays. Thus we get into a big range of
related questions and results on combinatorial complexity.
In \cite{aronovsharir} --- \cite{wiernik} one can
find useful ideas and
approaches to related questions.  %to such kind of problems we found in

%Most of result related to computational geometry, many contain asymptotical estimates.
%We give tight upper bounds for the number of vertices of a
%polygon that is the union or the intersection of two polygons with given numbers of
%convex and concave vertices. The similar question on graphs of the lower (or upper) envelope
%of two continuous piecewise linear functions considered.

For a piecewise-linear function  $f:\mathbb{R}\to \mathbb{R}$,
denote by $n(f)$ the total number of vertices in the graph of $f$,
denote by $c(f)$ the number of convex vertices.
%(here a vertex $a$ is said to be {\it convex} (or
%{\it concave}), if the graph of $f$ is a convex (concave) function in some
%neighborhood of $a$).
Thus $r(f)=n(f)-c(f)$ is the number of concave (reflex) vertices.

\begin{theor} \label{t1}
Let $n_i, c_i$, $i=1, 2, $ be given integers with $0\leq c_i\leq n_i$, $c_1\leq c_2$.
%For $i=1, 2$,
Consider  continuous piecewise-linear functions  $f_i:\mathbb{R}\to \mathbb{R}$ (for $i=1, 2$)
with $n(f_i)=n_i$, and  $c(f_i)=c_i$.
If $f_0$ is the lower envelope of $f_1$ and $f_2$, then
$n(f_0)\leq n_1+n_2+1+\min\{2c_1+n_2-c_2+1, c_1+c_2\}$, $c(f_0)\leq c_1+c_2$,
and these estimates are tight.
\end{theor}

\begin{corol} \label{t2}
Let $n_i$, $i=1, 2, $ be given non-negative integers.
%For $i=1, 2$,
Consider  continuous piecewise-linear functions $f_i:\mathbb{R}\to \mathbb{R}$ (for $i=1, 2$)
with $n(f_i)=n_i$.
If $f_0$  the lower envelope of
$f_1$ and $f_2$, then
$n(f_0)\leq 2n_1+2n_2+1- \left \lfloor \frac{|n_2-n_1|}{2} \right\rfloor$,
and this estimate is tight.
\end{corol}

Theorem \ref{t1} could be reformulated for the case of the upper envelope
by interchanging convex and concave vertices. It
could be iterated to obtain bounds for the number of
(convex, concave) vertices for lower envelopes of $k>2$ functions.

%Let $\{n_1, n_2, \ldots, n_k\}$ be a multiset of $k$ integers,
%$n_i\geq 3$ for $i\in \{1, 2, \ldots , k\}$.
%Consider a set $\{P_1, P_2, \ldots , P_k\}$ of be a set
%of polygons in the plane. Let $P_i$ be an $n_i$-gon $i\in \{1, 2, \ldots , k\}$.
%We set $S_i = P_i\cup Int (P_i)$.
%Let $S= \bigcap \limits_{i=1}^{k} S_i$
%Given that the boundary of  is an $n$-gon
%we estimate $n$.
%\vspace{1mm}

%Let $P$ be a polygon in the plane.
%We set $\overline{P}=P\cup Int(P)$.
%By $n(P)$ denote the number of vertices of $\overline{P}$,
%by $m(P)$ denote the number of interior angles of $\overline{P}$
%which are less then $\pi$.

Similar results (see below) obtained for unions and intersections of polygons in the plane.

For a polygon $P$, denote by $n(P)$  the total number of its vertices,
denote by $c(P)$ the number of its convex vertices.
Thus $r(P)=n(P)-c(P)$ is the number of concave vertices.

\begin{theor} \label{t1u}
Let $n_i, c_i$, $i=1, 2, $ be given integers with $3\leq c_i\leq n_i$, $c_1\leq c_2$.
%For $i=1, 2$,
Consider polygons $P_i$ (for $i=1, 2$)
with $n(P_i)=n_i$, and  $c(P_i)=c_i$.
If polygon $P_0$
is the union of $P_1$ and $P_2$, then
$n(P_0)\leq n_1+n_2+\min\{2c_1+n_2-c_2, c_1+c_2\}$, $c(P_0)\leq c_1+c_2$,
and these estimates are tight.
\end{theor}

\begin{corol} \label{t2u}
Let $n_i\geq 3$, $i=1, 2, $ be given integers.
%For $i=1, 2$,
Consider polygons (for $i=1, 2$)
with $n(P_i)=n_i$.
If polygon $P_0$
is the union of $P_1$ and $P_2$, then
$n(P_0)\leq 2n_1+2n_2- \left \lceil \frac{|n_2-n_1|}{2} \right\rceil$,
and this estimate is tight.
\end{corol}

In the following Theorem similar estimates presented for polygon $P_0$ that is
a connected component of $P_1\cap P_2$.
%the intersection
%$P_1\cap P_2$ or a subset of $P_1\cap P_2$ such that its boundary $\partial P_0$ is a subset
%$\partial P_1 \cup \partial P_2$.

\begin{theor} \label{t1p}
Let $n_i, r_i$, $i=1, 2, $ be given integers with $0\leq r_i\leq n_i-3$, $r_1\leq r_2$.
%For $i=1, 2$,
Consider polygons $P_i$ (for $i=1, 2$)
with $n(P_i)=n_i$, and  $r(P_i)=r_i$, $r_1\leq r_2$.
Let $P_0$ be a polygon such that $P_0 \subset P_1\cap P_2$
and $\partial P_0 \subset \partial P_1 \cup \partial P_2$.
We have $n(P_0)\leq n_1+n_2+\min\{2r_1+n_2-r_2, r_1+r_2\}$, $r(P_0)\leq r_1+r_2$,
and these estimates are tight.
\end{theor}

\begin{corol} \label{t2p}
Let $n_i\geq 3$, $i=1, 2, $ be given integers.
%For $i=1, 2$,
Consider polygons $P_i$ (for $i=1, 2$)
with $n(P_i)=n_i$.
Let $P_0$ be a polygon such that $P_0 \subset P_1\cap P_2$
and $\partial P_0 \subset \partial P_1 \cup \partial P_2$.
We have
$n(P_0)\leq 2n_1+2n_2-6-\max \left\{ \left\lfloor \frac{|n_1-n_2|}{2} \right\rfloor -1, 0\right\}$,
and this estimate is tight.
\end{corol}

%The formulated results could be iterated to obtain bounds for the number of
%(convex, concave) vertices for lower envelopes of $k>2$ functions.
%as well as for unions and intersections of $k>2$ polygons.

\smallskip

%\begin{theor} \label{t1p}
%Пусть $P_i$ --- плоский $n_i$-угольник, имеющий ровно $r_i$ углов, больших $\pi$, $i=1, 2$.
%Если пересечение $P_1\cap P_2$ является плоским $n$-угольником, то $n\leq n_1+n_2+r_1+r_2$.
%%причем $r\leq r_1+r_2$.
%\end{theor}

%\begin{theor} \label{t2p}
%Даны натуральные $n_1\geq 3$, $n_2\geq 3$.
%Рассматриваются плоские $n_i$-угольники $P_i$, $i=1, 2$,
%такие что пересечение $P_1\cap P_2$ является плоским $n$-угольником.
%Тогда наибольшее возможное значение $n$ равно $2n_1+2n_2-6-\max \left\{\left[\frac{|n_1-n_2|-2}{2}\right], 0\right\}$.
%%причем $r\leq r_1+r_2$.
%\end{theor}

{\bf Acknowledgement.}
%The author
%thanks the referee for a number of useful recommendations.
The author
thanks Prof. Herbert Edelsbrunner for pointing useful references
and Roman Karasev for helpful remarks and attention to this paper.

%Hi Roman,

%     if the assumption is that P cap Q is again a simple polygon
%   then the number of edges might be related to the Davenport-
%   Schinzel sequences ... and if I remember right, the order of
%   the seqence is such that the result is
%   a linear bound, perhaps twice the number of edges of P and Q,
%   or similar.

%     If I am on the right track, then your friend can find
%   pointers by following papers related to Davenport-Schinzel
%   sequences, mostly in the 80s and 90s by Sharir and coauthors.

%   Let me know if I was on the right track or the answer to the
%   problem is quite different from what I thought.

%   Greetings,
%     Herbert

%On Mon, 7 Nov 2011, Roman Karasev wrote:

%> Dear Herbert,
%>
%> This time I have a question on pure mathematics, not related to that grant
%> activity.
%>
%> One colleague of mine studied intersections of two polygons $P$ and $Q$ in
%> the plane. Polygons may be non-convex but supposed to be connected and simply
%> connected. Then he gives a tight upper bound on the number of edges of $P\cap
%> Q$ in terms of the number of edges of $P$ and $Q$ under the assumption that
%> $P\cap Q$ is again a polygon.
%>
%> What is known and what are the common references for such facts? Of course,
%> the assumption that $P\cap Q$ is connected and simply connected is not
%> natural for computational problems, but it general it sounds reasonable.
%>

\subsection{Proofs %of Theorems
}

\subsubsection{Proof of Theorem \ref{t1}}

Let $\Gamma _i$ be the graph of $f_i$, $i=0, 1, 2$.
If $\Gamma_1 \cap \Gamma_2 = \varnothing$,
then $\Gamma _0$ coincides to one of graphs  $\Gamma_1$, $\Gamma_2$; in this case
the statements of Theorem \ref{t1} is trivial.
Further suppose that $\Gamma_1 \cap \Gamma_2 \neq \varnothing$.

%%%%%define vertices

For $i=0, 1, 2$, by $C_i$, $R_i$  denote the sets of convex and concave vertices
in $\Gamma _i$, respectively. Thus $V_i = C_i\cup R_i$ is the set of all vertices in $\Gamma _i$.
We have $|V_i|=n_i$, $|C_i| = c_i$, $|R_i|=r_i=n_i-c_i$.
It is clear that $C_0 \subset C_1\cup C_2$, hence $c_0\leq c_1+c_2$.

%%%%%%define arcs

For points $a, b\in \Gamma _i$  ($i=0, 1, 2$) with $x$-coordinates $a_x<b_x$
by $[ab]_i$ denote the path
in $\Gamma _i$ connecting $a$ and $b$.
% such that
%the movement along $[xy]_i$ from $x$ to $y$ is counterclockwise ($\partial {P_i}$ is oriented
%counterclockwise).
Let $(ab)_i = [ab]_i\setminus \{a, b\}$. For
$a\in \Gamma_i$ define infinite paths $(a,+\infty)_i =
\{c \in \Gamma_i\, | \, c_x>a_x\}$ and $(-\infty,a)_i =
\{c \in \Gamma_i\, | \, c_x<a_x\}$.

%%%%%%%define $B$

Define an ordered set of {\it breakpoints} $B = \{ b_1, \ldots, b_k\}$
as a minimal subset of $\Gamma_1 \cap \Gamma_2$ satisfying the following
properties:

i) $x$-coordinates of $b_1, \ldots, b_k$ are ordered from the left to the right, i.~e.
$(b_1)_x < \ldots < (b_k)_x$;

ii) for each path $p\in \{(-\infty,b_1)_0, (b_1, b_2)_0, \ldots, (b_{k-1}, b_k)_0,
(b_k, +\infty)_0 \}$ either $p \subset \Gamma_1$
or $p \subset \Gamma_2$.

%Let us call $B$ the set of {\it the breakpoints }.
(In particular, if $\Gamma_1$ and $\Gamma_2$ intersect
{\it transversally} at each point of  $\Gamma_1 \cap \Gamma_2$,
then $B = \Gamma_1 \cap \Gamma_2$.)

%%%%%%%%%%%%

Let $K=\{1, 2, \ldots, k-1\}$.
Define subsets $K_1,  K_2\subset K$  as follows:
if $(b_tb_{t+1})_0 \subset \Gamma_1$, then we put $t\in K_1$, otherwise
(i.e. if $(b_tb_{t+1})_0 \subset \Gamma_2$) we put $t\in K_2$.

From the  condition of minimality in the definition of $B$ it easily follows that:

i) $K_1\cap K_2 =\varnothing $ (otherwise $(b_tb_{t+1})_0 \subset \Gamma_1 \cap \Gamma_2$
for some $t\in K$, hence $b_t$ could be removed from $B$);

ii) elements of $K_1$ and $K_2$ alter, i.~e.
$|K_1\cap \{t, t+1\}| = |K_2\cap \{t, t+1\}|=1$
for each $t\in K$
%of pairs $(1,2)$,
%$(2, 3)$, \ldots, $(k-1, k)$ contains exactly one number from $K_1$ and one number from $K_2$
(otherwise, if for instance, $\{t-1, t\} \subset K_1$, then
$(b_{t-1}, b_{t+1})_0 \subset  \Gamma_1$, and again $b_t$ could be removed from $B$).

From the last property we obtain

\begin{equation}  \label{alt1}
\left||K_1| - |K_2|\right|\leq 1.
\end{equation}

%\begin{lem} \label{topology}
%For $i=1, 2$ we have $\bigcup \limits _{t=1}^{k}(b_tb_{t+1})_i = \partial P_i \setminus B$.
%\end{lem}
%\dok
%\edok
%%%%%%%%%

%%%%%%%типы вершин...

%It is clear that $V_0\subset V_1\cup V_2\cup B$.
Define $K^c_1\subset K_1$ as follows:
$t\in K^c_1$ iff $t\in K_1$ and $C_1 \cap (b_t, b_{t+1})_0 \neq \varnothing$.
Take $K^r_1 = K_1\setminus K^c_1$.
Similarly, define $K^c_2\subset K_2$ and $K^r_2 =K_2\setminus K^c_2$. By definition,

\begin{equation}\label{K^c1}
|K^c_1|\leq c_1, \, \, \, |K^c_2|\leq c_2.
\end{equation}

\begin{lem} \label{reflex1}
If $t \in K^r_1$, then $(b_t b_{t+1})_2 \cap R_2\neq \varnothing$.
\end{lem}
\dok Let $t \in K^r_1$. This means that the path $[b_t b_{t+1}]_1=[b_t b_{t+1}]_0$
is the graph of a concave function.
% (that is the restriction of $f_1$ to $[(b_t)x (b_{t+1})x]$.
For all $x\in [(b_t)_x, (b_{t+1})_x]$ we have $f_2(x)\geq f_1(x)$,
$f_2((b_t)_x) =  f_1((b_t)_x)$, $f_2((b_{t+1})_x) =  f_1((b_{t+1})_x)$,
and there exists $z\in [(b_t)_x, (b_{t+1})_x]$ such that $f_2(z) > f_1(z)$.
This means that the path $[b_t b_{t+1}]_2$
%restriction of $f_2$ to
%$[(b_t)x (b_{t+1})x]$
is not a convex function, hence $(b_t b_{t+1})_2$
contains a concave vertex of $\Gamma_2$.
\edok

From Lemma \ref{reflex1} and the similar result for $K^r_2$ we obtain

\begin{equation}\label{K^r1}
|K^r_1| \leq r_2, \, \, \, |K^r_2| \leq r_1.
\end{equation}

From estimates (\ref{alt1}), (\ref{K^c1}), and (\ref{K^r1})
it follows that four
non-negative integers $x_i=|K^c_i|$, $y_i=|K^r_i|$, $i=1, 2$, satisfy the following
system:

\begin{equation}\label{cases1}
\begin{cases}
|(x_1+y_1)-(x_2+y_2)|\leq 1, \cr
x_1\leq c_1, x_2\leq c_2, \cr
y_1\leq r_2, y_2\leq r_1.
\end{cases}
\end{equation}

Note that $V_0 \subset V_1\cup V_2 \cup B$, and moreover,
$|V_0| \leq |V_1|+|V_2|+|B| -
\left| V_2 \cap \left(\bigcup _{t\in K_1} (b_t b_{t+1})_2\right) \right| -
\left| V_1 \cap \left(\bigcup _{t\in K_2} (b_t b_{t+1})_1\right) \right|$.
%From Lemma \ref{reflex1} it follows that for
%each $t \in K^r_1$ the path $(b_t b_{t+1})_2$ contains at least one vertex
%from $V_2\setminus V_0$. Similarly, for each $t \in K^r_2$ the path $(b_t b_{t+1})_1$
%contains at least one vertex from $V_1\setminus V_0$.
Thus Lemma \ref{reflex1} leads us to the following estimate:
$n_0 = |V_0| \leq |V_1|+|V_2|+|B|-|K^r_1|-|K^r_2|
\leq n_1+n_2+(x_1+y_1+x_2+y_2+1) - y_2 -y_1 =
n_1+n_2+x_1+x_2+1$. Hence $n_0\leq n_1+n_2+c_1+c_2+1$.
%If $2r_1+n_2-r_2+1 < r_1+r_2$, then the estimate could be improved
%in the following way:
From (\ref{cases1}) we have $x_2\leq x_1+y_1+1$, thus
$n_0 \leq n_1+n_2+x_1+x_2+1 \leq n_1+n_2+2x_1+y_1 +2  \leq
n_1+n_2+2c_1+r_2 +2 = n_1+n_2+2c_1+(n_2-c_2) +2$.
The estimates from the formulation of Theorem \ref{t1} established.

\vspace{3mm}

Now we show that the estimate is tight.
Given parameters $n_1, n_2$, and  $c_1\leq c_2$, we put
 $r_1=n_1-c_1, r_2=n_2-c_2 $ and
construct the corresponding example  $E(c_1, c_2, r_1, r_2)$ that is the required configuration
of graphs $(\Gamma_1, \Gamma_2)$. In the examples below we have $C_0=C_1\cup C_2$, $C_1\cap C_2 = \varnothing$,
and hence $c_0=c_1+c_2$.

%Let $A$ be a finite set of points, $|A|\geq 2$.
Let $A = \{a_0, a_1, \ldots, a_t\}$ be a finite set of points, $t\geq 2$, and $a_i$ are ordered by $x$-coordinate, i.~e.
$(a_0)_x < (a_1)_x < \ldots < (a_t)_x$.
By $\gamma (A)$ denote the union of segments $[a_0, a_1], \ldots, [a_{t-1}, a_t]$.
By $\Gamma (A)$ denote the union of two rays
$[a_1 a_0)$, $[a_{t-1} a_t)$ and $\gamma (A)$.

We start with $E(c_1, c_2, 0, 0)$.
On a graph of a strictly convex function $y=g(x)$ take points $a^i_0, a^i_1, \ldots,
a^i_{c_i+1}$, $i=1, 2$, such that their $x$-coordinates satisfy
inequalities $(a^2_0)_x< (a^1_0)_x < (a^2_1)_x < (a^1_1)_x < \ldots < (a^2_{c_1+1})_x <
(a^1_{c_1+1})_x$,
and (for $c_2>c_1$) $(a^2_{c_1+1})_x < (a^2_{c_1+2})_x < \ldots < (a^2_{c_2})_x
< (a^1_{c_1+1})_x < (a^2_{c_2 +1})_x$.
Define $\Gamma _i = \Gamma (A^i)$, where $A^i = \{ a^i_0, \ldots, a^i_{c_i+1}\}$,
$i=1, 2$.
In this example $n_0=n_1+n_2+k$, where $k$ is the number of breakpoints.
We have $2c_1+1$ breakpoints $b_1, b_2, \ldots, b_{2c_1+1}$
with $(a^1_0)_x < (b_1)_x< (a^2_1)_x < (b_2)_x < (a^1_1)_x <  (b_3)_x  <
\ldots < (b_{2c_1+1})_x < (a^2_{c_1+1})_x $. If $c_1=c_2$, then this is the complete list of
breakpoints.
If $c_1<c_2$, then we have one more breakpoint $b_{2c_1+2}$  with
$(a^2_{c_2})_x < (b_{2c_1+2})_x< (a^1_{c_1+1})_x $. Thus $k=\min \{2c_1+2, c_1+c_2+1\}$, as
required.

Now to obtain
 $E(c_1, c_2, 0, r_2)$, where $0\leq r_2\leq c_2-c_1-1$,
 let us modify $\Gamma _2$ in $E(c_1, c_2, 0, 0)$.
Define points $d_i$, $i=1, 2, \ldots, r_2$, such that
$(a^2_{c_1+i})_x <(d_i)_x < (a^2_{c_1+i+1})_x$ and $(d_i)_y> f_1((d_i)_x)$.
Replace $\Gamma _2 = \Gamma (A^2)$ by $\Gamma (A^2\cup \{d_1, \ldots, d_{r_2}\})$.
%the segment $[a^2_{c_1+i}a^2_{c_1+i+1}]$ by the union of segments
%$[a^2_{c_1+i}d_i]\cup[d_ia^2_{c_1+i+1}]$, for $i=1, 2, \ldots, r_2$.
 As a result from $E(c_1, c_2, 0, 0)$ we obtain $E(c_1, c_2, 0, r_2)$ with $n_0$
increased by $2r_2$, as required ($2r_2$ additional breakpoints appeared).

Finally, to construct $E(c_1, c_2, r_1, r_2)$ for any $(c_1, c_2, r_1, r_2)$
it is sufficient to
take $E(c_1, c_2, 0, z)$ constructed above, where $z=\min \{r_2, \max\{0, c_2-c_1-1\}\}$, and
modify $\Gamma _1$ and $\Gamma _2$ by adding  $r_1$ and $r_2-z$ concave vertices, respectively,
to get $n_0$ increased by $r_1 +r_2-z$. Adding a concave vertex to $\Gamma_1$ (for
$\Gamma_2 $ the procedure is analogous) could be done in the following way:
take a point $e $ on $\gamma_1\cap \Gamma_0 $,
slightly move $e$ up to a new position $e'$, and replace $\Gamma_1 = \Gamma (A)$
by $\Gamma (A\cup \{e'\})$.

\subsubsection{Proof of Corollary \ref{t2}}

Now $n_1$, $n_2$, are fixed, while $c_1$ and $c_2$  could vary.
Suppose $n_1\leq n_2$, and put $m=n_2-n_1$.

For fixed $c_1, c_2$ from Theorem \ref{t1} we know the tight
upper  bound $N(n_1, n_2, c_1, c_2)$ for $n(f_0)$. Thus it remains just to find the
maximum of $N(n_1, n_2, x, y)$ while $(x, y)$ runs over $\{0, 1, \ldots, n_1\}\times
\{0, 1, \ldots, n_2\}$.

If $x\geq y$, then
$N(n_1, n_2, x, y)  = n_1+n_2+1+\min\{2y+n_1-x+1, x+y\} \leq
n_1+n_2+1+\min\{2y+n_2-x+1, x+y\} = N(n_1, n_2, y, x)$. Thus further we
assume that $x\leq y$.

If $x<y$ and $x<n_1$, then
$N(n_1, n_2, x, y)  = n_1+n_2+1+\min\{2x+n_2-y+1, x+y\} < N(n_1, n_2, x+1, y)$.
Thus it remains to consider cases $x=y$ and $x=n_1$.

Since $N(n_1, n_2, x, x)  = n_1+n_2+1+\min\{x+n_2+1, 2x\} \leq N(n_1, n_2, n_1, n_1)$,
it remains to consider $N(n_1, n_2, n_1, y)$, $n_1\leq y\leq n_2$.

%We have $N(n_1, n_2, x, x)  = n_1+n_2+1+\min\{x+n_2+1, 2x\} \leq n_1+n_2+1+2x
%\leq 3n_1+n_2+1  = 2n_1+2n_2+1 -m \leq 2n_1+2n_2+2- \left \lceil \frac{m+1}{2} \right\rceil$.

Let $y=n_1+z$, $z\in \{0, \ldots, m\}$.
We have $N(n_1, n_2, n_1, y)  = n_1+n_2+1+\min\{2n_1+n_2-n_1-z+1, 2n_1+z\} =
3n_1+n_2+1 + \min\{m+1-z, z\} \leq 3n_1+n_2+1 + \left\lfloor \frac{m+1}{2} \right\rfloor
= 2n_1+2n_2+1 -m + \left\lfloor \frac{m+1}{2} \right\rfloor =
 %= 2n_1+2n_2+2 - \left \lceil \frac{m+1}{2} \right\rceil
  2n_1+2n_2+1 - \left \lfloor \frac{m}{2} \right\rfloor$.
% $ = 2n_1+2n_2+2- \left \lceil \frac{|n_2-n_1|+1}{2} \right\rceil$.
The last inequality
turns to equality for $x=n_1$ and $y = n_1+ \left\lfloor \frac{m+1}{2} \right\rfloor$.

\subsubsection{Proof of Theorem~\ref{t1u}}

Mainly we follow  the outline of the proof of Theorem~\ref{t1}.

%If $\partial {P_1} \cap \partial {P_2} = \varnothing$,
%then $P_0$ coincides to one of polygons $P_1$, $P_2$; in this case
%the statements of Theorem~\ref{t1u} is trivial.
%Further suppose that $\partial {P_1} \cap \partial {P_2} \neq \varnothing$.

%%%%%define vertices

Let $P_i$ ($i=0, 1, 2$) be polygons such that $P_0=P_1\cup P_2$.
For $i=0, 1, 2$, by $C_i$, $R_i$  the set of convex and concave %(reflex)
vertices
in $P_i$, respectively. Thus $V_i = C_i\cup R_i$ is the set of all vertices in $P_i$.
We have $|V_i|=n_i$, $|C_i| = c_i$, $|R_i|=r_i=n_i-c_i$.
It is clear that $C_0 \subset C_1\cup C_2$, hence $c_0\leq c_1+c_2$.

%%%%%%define arcs

For points $a, b\in \partial {P_i}$ ($i=0, 1, 2$), by $[ab]_i$ denote the path
in $\partial {P_i}$ connecting $a$ and $b$ such that
the movement along $[ab]_i$ from $a$ to $b$ is counterclockwise ($\partial {P_i}$ is supposed
to be oriented). We set $(ab)_i = [ab]_i\setminus \{a, b\}$.

%%%%%%%define $B$

Define an ordered set of breakpoints  $B = \{ b_1, b_2, \ldots, b_k\}$
as a minimal subset of $\partial {P_1} \cap \partial {P_2}$ satisfying the following
conditions:

i) if $t, s\in \{1, \ldots, k\}$, $t\neq s$, then
$(b_tb_{t+1})_0 \cap (b_sb_{s+1})_0 = \varnothing $, $i=1, 2$
%$\bigcup \limits _{t=1}^{k}(b_tb_{t+1})_0 = \partial P_0 \setminus B$
(here and further we set $b_{k+1}=b_1$),
in other words, $b_1, b_2, \ldots, b_k$ are arranged on $\partial P_0$ in
the counterclockwise order.

ii) for each $t\in \{1, 2, \ldots, k\}$ we have either $(b_tb_{t+1})_0 \subset \partial P_1$
or $(b_tb_{t+1})_0 \subset \partial P_2$.

%(Let us call $B$ the set of {\it the breakpoints }. If $\partial P_1$ and $\partial P_2$ intersect
%transversally at each point, then $B = \partial {P_1} \cap \partial {P_2}$.)

%%%%%%%%%%%%

Further we assume that $k=|B|\geq 2$, otherwise we have either $\partial P_0 = \partial P_1$ or
$\partial P_0 = \partial P_2$, and the statement of Theorem is trivial.

Let $K=\{1, 2, \ldots, k\}$.
Define subsets $K_1,  K_2\subset K$  as follows:
if $(b_tb_{t+1})_0 \subset \partial P_1$, then we put $t\in K_1$, otherwise
(i.e. if $(b_tb_{t+1})_0 \subset \partial P_2$) we put $t\in K_2$.

From the  condition of minimality in definition of $B$ it easily follows that:

i) $K_1\cap K_2 =\varnothing $;
%(otherwise $(b_tb_{t+1})_0 \subset \partial P_1 \cap \partial P_2$
%for some $t\in K$, hence $b_t$ could be removed from $B$);

ii) elements of $K_1$ and $K_2$ alter, i.~e.
$|K_1\cap \{t, t+1\}| = |K_2\cap \{t, t+1\}|=1$
for each $t\in K$.
%
%i.~e.  each of pairs $(1,2)$,
%$(2, 3)$, \ldots, $(k, 1)$ contains exactly one number from $K_1$ and one number from $K_2$
%(otherwise, if for instance, $t-1, t \in K_1$ then
%$(b_{t-1}, b_{t+1})_0 \subset  \partial P_1$, and again $b_t$ could be removed from $B$).
From the last property we obtain
\begin{equation}  \label{alt2}
|K_1| = |K_2|,
\end{equation}
in particular, $k$ is even.

\begin{lem} \label{topology}
If $t, s\in K$, $t\neq s$, then
$(b_tb_{t+1})_i \cap (b_sb_{s+1})_i = \varnothing $, $i=1, 2$.
\end{lem}
\dok %For $i=1, 2$ we have $\bigcup \limits _{t=1}^{k}(b_tb_{t+1})_i = \partial P_i \setminus B$....
Let $i\in \{1, 2\}$, and $t\in K$. We need to prove that
$(b_tb_{t+1})_i\cap B= \varnothing$. If $t\in K_i$, then this is obvious, since
$(b_tb_{t+1})_i = (b_tb_{t+1})_0$. If $t\notin K_i$, then $(b_{t-1}b_{t})_i = (b_{t-1}b_{t})_0$.
Suppose to the contrary that $b_s\in (b_tb_{t+1})_i$. Then $b_{t+1}\notin (b_tb_s)_i$.
The path $(b_tb_s)_i$ divides $P_0$ into polygons; for one of these polygons $P'$ we have
$\partial P' = (b_tb_s)_i \cup (b_sb_t)_0$.
Thus $(b_{t-1}b_{t})_i \subset \partial P'$, hence $P_i\subset P'$; but
$b_{t+1}  \in \partial P_0$, $b_{t+1}  \notin \partial P'$, hence
$b_{t+1}  \notin P'$. This is the contradiction.
\edok
%%%%%%%%%

%%%%%%%типы вершин...

%It is clear that $V_0\subset V_1\cup V_2\cup B$.
Further define the subset $K^c_1\subset K_1$ as follows:
$t\in K^c_1$ iff $t\in K_1$ and $C_1 \cap (b_t, b_{t+1})_1 \neq \varnothing$.
Take $K^r_1 = K_1\setminus K^c_1$.
Similarly, define $K^c_2\subset K_2$ and $K^r_2 =K_2\setminus K^c_2$. By definition,

\begin{equation}\label{K^c}
|K^c_1|\leq c_1, \, \, \, |K^c_2|\leq c_2.
\end{equation}

\begin{lem} \label{reflex}
If $t \in K^r_1$, then $(b_t b_{t+1})_2 \cap R_2\neq \varnothing$.
\end{lem}
\dok Let $t \in K^r_1$. This means that the path $(b_t b_{t+1})_1=(b_t b_{t+1})_0$
could contain only concave vertices of $P_1$.
Consider points $x \in [b_t b_{t+1})_0\cap \partial P_2$ and
$y \in (x b_{t+1}]_0 \cap \partial P_2$
such that $(xy)_1\cap (xy)_2 = \varnothing$. Consider the polygon $M$ such that
$\partial M = [xy]_1\cup [xy]_2$. Note that $(xy)_1$ contains no convex vertices of $M$.
From the other hand, $M$ contains at least 3 convex vertices, hence
%there exist three convex vertices in $[xy]_2$, and hence
$(xy)_2 $ contains at least one convex vertex $z$.
This vertex $z\in (xy)_2 \subset (b_t b_{t+1})_2$ is a concave vertex of $P_2$.
\edok

From Lemma \ref{reflex} and the similar result for $K^r_2$ using Lemma \ref{topology}
 we obtain
\begin{equation}\label{K^r}
|K^r_1| \leq r_2, \, \, \, |K^r_2| \leq r_1.
\end{equation}

From (\ref{alt2}), (\ref{K^c}), and (\ref{K^r})
it follows that four
non-negative integers $x_i=|K^c_i|$, $y_i=|K^r_i|$, $i=1, 2$, satisfy the following
system:

\begin{equation}\label{cases2}
\begin{cases}
x_1+y_1=x_2+y_2, \cr
x_1\leq c_1, x_2\leq c_2, \cr
y_1\leq r_2, y_2\leq r_1.
\end{cases}
\end{equation}

Note that $V_0 \subset V_1\cup V_2 \cup B$, and moreover,
$|V_0| \leq |V_1|+|V_2|+|B| -
\left| V_2 \cap \left(\bigcup_{t\in K_1} (b_t b_{t+1})_2)\right) \right| -
\left| V_1 \cap \left(\bigcup_{t\in K_2} (b_t b_{t+1})_1\right) \right|$.
%From Lemma \ref{reflex} it follows that for
%each $t \in K^r_1$ the path $(b_t b_{t+1})_2$ contains at least one vertex
%from $V_2\setminus V_0$. Similarly, for each $t \in K^r_2$ the path $(b_t b_{t+1})_1$
%contains at least one vertex from $V_1\setminus V_0$.
Thus Lemma \ref{reflex} leads us to the following estimate:
$n_0 = |V_0| \leq |V_1|+|V_2|+|B|-|K^r_1|-|K^r_2|
\leq n_1+n_2+(x_1+y_1+x_2+y_2) - y_2 -y_1 =
n_1+n_2+x_1+x_2$. Hence $n_0\leq n_1+n_2+c_1+c_2$.
%If $2r_1+n_2-r_2+1 < r_1+r_2$, then the estimate could be improved
%in the following way:
From (\ref{cases2}) we obtain $x_2\leq x_1+y_1$, thus we have
$n_0 \leq n_1+n_2+x_1+x_2 \leq n_1+n_2+2x_1+y_1  \leq
n_1+n_2+2c_1+r_2 = n_1+n_2+2c_1+(n_2-c_2)$.
The estimates from the formulation of Theorem \ref{t1u} established.

\vspace{3mm}

Now we show that the estimate is tight.
Given parameters $n_1, n_2$, and  $c_1\leq c_2$, we put
 $r_1=n_1-c_1, r_2=n_2-c_2 $ and
construct the corresponding example  $U(c_1, c_2, r_1, r_2)$ that is the required configuration
of polygons $P_1$ and $P_2$. In the examples below we have $C_0=C_1\cup C_2$, $C_1\cap C_2 = \varnothing$,
and hence $c_0=c_1+c_2$.

%Let $A$ be a finite set of points, $|A|\geq 2$.
%Let $A = \{a_0, a_1, \ldots, a_t\}$, where $a_i$ are ordered by $x$-coordinate, i.~e.
%$(a_0)_x < (a_1)_x < \ldots < (a_t)_x$. By $\Gamma (A)$ denote
%the graph that is the union of two rays
%$[a_1 a_0)$, $[a_{t-1} a_t)$ and segments $[a_0, a_1], \ldots, [a_{t-1}, a_t]$ .

We start with $U(c_1, c_2, 0, 0)$.
To define polygons $P_i = a^i_1 \ldots a^i_{c_i}$, $i=1, 2$,
 take points $a^i_1, \ldots,
a^i_{c_i}$ on a  circle in the counterclockwise cyclic order
 $a^2_1, a^1_1 , a^2_2, a^1_2 ,\ldots , a^2_{c_1} , a^1_{c_1} ,
 a^2_{c_1+1} , \ldots,  a^2_{c_2} ,  a^2_1$.
In this example $n_0=n_1+n_2+k$, where $k=2c_1$ is the required number of breakpoints.

In the construction $U(c_1, c_2, 0, 0)$
take a point $q$ in the interior of
 $P_1\cap P_2$. Modify $P _2$ to obtain
 $U(c_1, c_2, 0, r_2)$, where $r_2\leq c_2-c_1$ in the following way.
Define points $d_i$ in the interior of the segment $[qa^2_{c_1+i}]$, $i=1, 2, \ldots, r_2$,
such that each segment $[d_ia^2_{c_1+i}]$ intersects the segment $[a^1_{c_1}a^1_1]$.
In $\partial P_2$ replace
the segment $[a^2_{c_1+i}a^2_{c_1+i+1}]$ by the union of segments
$[a^2_{c_1+i}d_i]\cup[d_ia^2_{c_1+i+1}]$, for $i=1, 2, \ldots, r_2$ (here $a^2_{c_2+1} = a^2_{1}$).
 As a result we obtain $E(c_1, c_2, 0, r_2)$ with $n_0$
increased by $2r_2$.

Finally, to construct $U(c_1, c_2, r_1, r_2)$ for any $(c_1, c_2, r_1, r_2)$
it is sufficient to
take $U(c_1, c_2, 0, z)$ constructed above, where $z=\min \{r_2, c_2-c_1\}$,  and
modify $P _1$ and $P _2$ by adding  $c_1$ and $c_2-z$ concave vertices, respectively,
to get $n_0$ increased by $r_1 +r_2-z$. Adding a concave vertex to $P_1$
 could be done in the following way.
Take a point $e$ of the segment $[ab]\subset \partial P_1$, where $a, b \in V_1$,
 $e\notin  P_2 $.
In $\partial P_1$ replace $[ab]$ by $[ae']\cup [e'b]$, where $e'$ obtained from $e$
by a slight shift inside $P_1$.

\subsubsection{Proof of Corollary \ref{t2u}}

Now $n_1$, $n_2$ are fixed, while $c_1$ and $c_2$  could vary.
Suppose $n_1\leq n_2$, and put $m=n_2-n_1$

For fixed $c_1, c_2$, from Theorem \ref{t1u} we know the tight
upper  bound $N'(n_1, n_2, c_1, c_2)$ for $n(P_0)$, thus it remains to find the
maximum of $N'(n_1, n_2, x, y)$ while $(x, y)$ runs over $\{3, \ldots, n_1\}\times
\{3, \ldots, n_2\}$.

If $x\geq y$, then
$N'(n_1, n_2, x, y)  = n_1+n_2+\min\{2y+n_1-x, x+y\} \leq
n_1+n_2+\min\{2y+n_2-x, x+y\} = N'(n_1, n_2, y, x)$. Thus further we
assume that $x\leq y$.

If $x<y$ and $x<n_1$, then
$N'(n_1, n_2, x, y)  = n_1+n_2+\min\{2x+n_2-y, x+y\} < N'(n_1, n_2, x+1, y)$.
Thus it remains to consider cases $x=y$ and $x=n_1$.

We have $N'(n_1, n_2, x, x)  = n_1+n_2+\min\{x+n_2, 2x\} \leq N'(n_1, n_2, n_1, n_1).$
It remains to consider $N'(n_1, n_2, n_1, y)$, $n_1\leq y\leq n_2$.

Let $y=n_1+z$, $z\in \{0, \ldots, m\}$.
We have $N'(n_1, n_2, n_1, y)  = n_1+n_2+\min\{2n_1+n_2-n_1-z, 2n_1+z\} =
3n_1+n_2+\min\{m-z, z\} \leq 3n_1+n_2+\left\lfloor \frac{m}{2} \right\rfloor
= 2n_1+2n_2-m + \left\lfloor \frac{m}{2} \right\rfloor =
 2n_1+2n_2- \left \lceil \frac{m}{2} \right\rceil$.
% $ = 2n_1+2n_2+2- \left \lceil \frac{|n_2-n_1|+1}{2} \right\rceil$.
The last inequality
turns to equality for $x=n_1$ and $y = n_1+ \left\lfloor \frac{m}{2} \right\rfloor$.

\subsubsection{Proof of Theorem \ref{t1p}}

The situation here is close  to Theorem~\ref{t1u} (concepts of interior
and exterior interchange, the same for concave and convex vertices).

Let $P_0$ be a polygon such that $P_0 \subset P_1\cap P_2$
and $\partial P_0 \subset \partial P_1 \cup \partial P_2$.
We keep all the the other notation from the proof of Theorem \ref{t1u} up to Lemma \ref{topology}
(i.~e. define an ordered set of breakpoints  $B = \{ b_1, b_2, \ldots, b_k\}$,
define subsets $K_1$, $K_2$ of $K=\{1, \ldots, k\}$).
Now we have $R_0\subset R_1\cup R_2$, and hence $r_0\leq r_1+r_2$.
Again we have

\begin{equation}  \label{alt3}
|K_1| = |K_2|.
\end{equation}

The following Lemmas analogous to Lemmas \ref{topology} and \ref{reflex}.

\begin{lem} \label{topology2}
If $t, s\in K$, $t\neq s$, then
$(b_tb_{t+1})_i \cap (b_sb_{s+1})_i = \varnothing $, $i=1, 2$.
\end{lem}
\dok
Let $i\in \{1, 2\}$, and $t\in K$. We need to prove that
$(b_tb_{t+1})_i\cap B= \varnothing$. If $t\in K_i$, then this is obvious, since
$(b_tb_{t+1})_i = (b_tb_{t+1})_0$. If $t\notin K_i$, then $(b_{t-1}b_{t})_i = (b_{t-1}b_{t})_0$.
Suppose to the contrary that $b_s\in (b_tb_{t+1})_i$. Then $b_{t+1}\notin (b_tb_s)_i$.
%The path $(b_tb_s)_i$ divides $P_0$ into polygons; for one of these polygons
Consider the polygon $P'$ such that
$\partial P' = (b_tb_s)_i \cup (b_sb_t)_0$.
We have $P_0\subset P' \subset P_i$, and
$b_{t+1}  \notin \partial P'$. Hence $b_{t+1}$ is the interior point of $P_i$.
This is the contradiction.
\edok

%%%%%%%типы вершин...

%It is clear that $V_0\subset V_1\cup V_2\cup B$.
Now define $K^r_1\subset K_1$ as follows:
$t\in K^r_1$ iff $t\in K_1$ and $R_1 \cap (b_t, b_{t+1})_1 \neq \varnothing$.
Take $K^c_1 = K_1\setminus K^r_1$.
Similarly, define $K^r_2\subset K_2$ and $K^c_2 =K_2\setminus K^r_2$. By definition,

\begin{equation}\label{K^r2}
|K^r_1|\leq r_1, \, \, \, |K^r_2|\leq r_2.
\end{equation}

\begin{lem} \label{reflex2}
If $t \in K^c_1$, then $(b_t b_{t+1})_2 \cap C_2\neq \varnothing$.
\end{lem}
\dok Let $t \in K^c_1$. This means that the path $(b_t b_{t+1})_1=(b_t b_{t+1})_0$
could contain only convex vertices of $P_1$.
Consider points $x \in [b_t b_{t+1})_0\cap \partial P_2$ and
$y \in (x b_{t+1}]_0 \cap \partial P_2$
such that $(xy)_1\cap (xy)_2 = \varnothing$. Consider the polygon $M$ with
$\partial M = [xy]_1\cup [xy]_2$. Note that $(xy)_1$ contains no convex vertices of $M$.
From the other hand, $M$ contains at least 3 convex vertices, hence
%there exist three convex vertices in $[xy]_2$, and hence
$(xy)_2 $ contains at least one convex vertex $z$.
This vertex $z\in (xy)_2 \subset (b_t b_{t+1})_2$ is a convex vertex of $P_2$.
\edok

From Lemma \ref{reflex2} and the similar result for $K^r_2$ using Lemma \ref{topology2}
 we obtain
\begin{equation}\label{K^c2}
|K^c_1| \leq c_2, \, \, \, |K^c_2| \leq c_1.
\end{equation}

From (\ref{alt3}), (\ref{K^r2}), and (\ref{K^c2})
it follows that four
non-negative integers $x_i=|K^r_i|$, $y_i=|K^c_i|$, $i=1, 2$, satisfy the following
system:

\begin{equation}\label{cases3}
\begin{cases}
x_1+y_1=x_2+y_2, \cr
x_1\leq r_1, x_2\leq r_2, \cr
y_1\leq c_2, y_2\leq c_1.
\end{cases}
\end{equation}

Note that $V_0 \subset V_1\cup V_2 \cup B$, and moreover,
$|V_0| \leq |V_1|+|V_2|+|B| -
\left| V_2 \cap \left(\bigcup_{t\in K_1} (b_t b_{t+1})_2)\right) \right| -
\left| V_1 \cap \left(\bigcup_{t\in K_2} (b_t b_{t+1})_1\right) \right|$.
%From Lemma \ref{reflex} it follows that for
%each $t \in K^c_1$ the path $(b_t b_{t+1})_2$ contains at least one vertex
%from $V_2\setminus V_0$. Similarly, for each $t \in K^c_2$ the path $(b_t b_{t+1})_1$
%contains at least one vertex from $V_1\setminus V_0$.
Thus Lemma \ref{reflex2} leads us leads us to the following estimate:
$n_0 = |V_0| \leq |V_1|+|V_2|+|B|-|K^c_1|-|K^c_2|
\leq n_1+n_2+(x_1+y_1+x_2+y_2) - y_2 -y_1 =
n_1+n_2+x_1+x_2$. Hence $n_0\leq n_1+n_2+r_1+r_2$.
%If $2r_1+n_2-r_2+1 < r_1+r_2$, then the estimate could be improved
%in the following way:
From (\ref{cases3}) we have $x_2\leq x_1+y_1$, hence
$n_0 \leq n_1+n_2+x_1+x_2 \leq n_1+n_2+2x_1+y_1  \leq
n_1+n_2+2r_1+c_2 = n_1+n_2+2r_1+(n_2-c_2)$.
The estimates from the formulation of Theorem \ref{t1p} established.

\vspace{3mm}

Now we show that the estimate is tight.
Given parameters $n_1, n_2$, and  $r_1\leq r_2$, we put
 $c_1=n_1-r_1, c_2=n_2-r_2 $, and
construct the required configuration
of polygons $(P_1, P_2)$ such that $P_0=P_1\cap P_2$ is a polygon with
$n_0=M(r_1, r_2, c_1, c_2)$, where
$M(r_1, r_2, c_1, c_2)=r_1+r_2+c_1+c_2+\min\{2r_1+c_2, r_1+r_2\}$.
In the examples below we have $R_0=R_1\cup R_2$, $R_1\cap R_2 = \varnothing$,
and hence $r_0=r_1+r_2$.

%Let $A$ be a finite set of points, $|A|\geq 2$.
%Let $A = \{a_0, a_1, \ldots, a_t\}$, where $a_i$ are ordered by $x$-coordinate, i.~e.
%$(a_0)_x < (a_1)_x < \ldots < (a_t)_x$. By $\Gamma (A)$ denote
%the graph that is the union of two rays
%$[a_1 a_0)$, $[a_{t-1} a_t)$ and segments $[a_0, a_1], \ldots, [a_{t-1}, a_t]$ .

We use the notation from the construction of examples  in the proof of
Theorem \ref{t1}. We start with the configuration $E(c_1', c_2', r_1', r_2')$ of graphs
$(\Gamma_1, \Gamma_2 )$ having lower envelope with
$M(c_1', c_2', r_1', r_2')=c_1'+c_2'+r_1'+r_2'+1+\min\{2c_1'+r_2'+1, c_1'+c_2'\}$ vertices,
where we set $r_1'=c_1-3$, $r_2'=c_2-3$, $c_1'=r_1$, and $c_2'=\begin{cases}
r_2, \,\, \, \, \, \, \, \, \, \, \, \, \mbox{if} \, \, r_2\leq r_1+1, \cr
r_2-1, \, \, \mbox{if} \, \,  r_2= r_1+2, \cr
r_2-2, \, \, \mbox{if} \, \, r_2\geq r_1+3.
\end{cases}$
By construction, $\Gamma _i = \Gamma (A'^i)$, where $A'^i$ is a finite subset of
$\{(x,y)\, | \, y\geq g(x) \}$ (here $g(x)$ is a strictly convex function).
%points $a$ with $a_y\geq g(a_x)$.
Recall that $a^i_0 $ and $a^i_{c'_i+1}$ are points of $A'^i$
with minimal and maximal $x$-coordinate, respectively; to simplify further notation we set
$a^i= a^i_0$ and $z^i=a^i_{c'_i+1}$.

Take a point $q$ with sufficiently large negative $y$-coordinate
($q_y<<0$) and define $P_1$ as a polygon with
$\partial P_1 = \gamma (A'^1) \cup [a^1q] \cup [z^1q] $.
Note that $r(P_1) = c_1'$, $c(P_1)=r_1'+3$, since convex vertices of $\Gamma_1$ turn to
concave vertices of $P_1$ (and vice versa), and 3 more convex vertices
$a^1$, $z^1$, and $q$ appeared.

In the interior of $P_1$ choose a point $p$ with $p_y<<0$ (choose $p$ close to $q$) so that
$\gamma (A'^2) \cup [a^2p] \cup [z^2p] $ is a boundary of a polygon.
In the segments $[a^2p]$ and $[z^2p]$ take point points $p_1$ and $p_2$, respectively.
Suppose that $p_1$ and $p_2$ are placed close to $p$ so that $p_1, p_2$ belong to the interior of
$P_1$. Take point $p_0$ so that $q$ belongs to the interior of the segment $[pp_0]$. Let
$p_3$ be the common point of lines $a^2p$ and $z^2p_0$. To define $P_2$ consider four cases.

For $r_2\geq r_1+3 $, define $P_2$ as a polygon with
$\partial P_2 = \gamma (A'^2) \cup [a^2p_1] \cup [p_1p_0] \cup [p_0p_2] \cup [p_2z^2] $.
Note that $r(P_2) = c_2'+2=r_2$, $c(P_2)=c_1'+3$, since convex vertices of $\Gamma_2$ turn to
concave vertices of $P_1$ (and vice versa),  3 more convex vertices
$a^2$, $z^2$, $p_0$ and 2 more concave vertices
$p_1$, $p_2$ appeared.
In this construction $V_0$ contains all
vertices of $\Gamma_0 = \Gamma_1\cap \Gamma_2$ and 9 additional vertices
$a^1, p_1, q, p_2, z^1, [a^2p_1]\cap [a^1q], [p_1p_0]\cap [a^1q],
[p_0p_2]\cap [z^1q], [z^2p_2]\cap [z^1q]$, thus
$|V_0|= r_1+(r_2-2)+(c_1-3)+(c_2-3)+1+
\min \{2r_1+(c_2-3)+1, r_1+(r_2-2)\}+9 = M(r_1, r_2, c_1, c_2)$.

For $r_2= r_1+2 $, define $P_2$ as a polygon with
$\partial P_2 = \gamma (A'^2) \cup [a^2p_3] \cup [p_3p_2] \cup [p_2z^2] $.
Note that $r(P_2) = c_2'+1=r_2$, $c(P_2)=c_1'+3$ (3 convex vertices
$a^2$, $z^2$, $p_3$ and 1 concave vertex $p_2$ appeared).
In this case $V_0$ contains all
vertices of $\Gamma_0$ and 7 additional vertices
$a^1, p_2, a^1, [a^2p_3]\cap [a^1q], [a^2p_3]\cap [z^1q],
[p_3p_2]\cap [z^1q], [z^2p_2]\cap [z^1q]$, thus
$|V_0|= r_1+(r_2-1)+(c_1-3)+(c_2-3)+1+
(r_1+(r_2-1))+7 = M(r_1, r_2, c_1, c_2)$.

For $r_2= r_1+1 $, define $P_2$ as a polygon with
$\partial P_2 = \gamma (A'^2) \cup [a^2p] \cup [pz^2] $.
Here $r(P_2) = c_2'=r_2$, $c(P_2)=c_1'+3$ (3 convex vertices
$a^2$, $z^2$, $p$ appeared).
In this case $V_0$ contains all vertices of $\Gamma_0$ and 5 additional vertices
$a^1, p, z^1, [a^2p]\cap [a^1q], [z^2p]\cap [z^1q]$, thus
$|V_0|= r_1+r_2+(c_1-3)+(c_2-3)+1+
(r_1+r_2)+5 = M(r_1, r_2, c_1, c_2)$.

Finally, For $r_2= r_1$, define $P_2$ as a polygon with
$\partial P_2 = \gamma (A'^2) \cup [a^2p_3] \cup [p_3z^2] $.
Here $r(P_2) = c_2'=r_2$, $c(P_2)=c_1'+3$ (3 convex vertices
$a^2$, $z^2$, $p_3$ appeared).
In this case $V_0$ contains all vertices of $\Gamma_0$ and 5 additional vertices
$a^1, z^2, [a^2p_3]\cap [a^1q],
[a^2p_3]\cap [z^1q], [z^2p_3]\cap [z^1q]$
(note that in this case $z^2$ belongs to the interior of $P_1$), thus
$|V_0|=  r_1+r_2+(c_1-3)+(c_2-3)+1+
(r_1+r_2)+5 = M(r_1, r_2, c_1, c_2)$.

\subsubsection{Proof of Corollary \ref{t2p}}

Now we need to maximize $N'(n_1, n_2, x, y)$
for fixed $n_1$, $n_2$, while $(r_1, r_2)$  vary over
$\{0,  \ldots, n_1-3\}\times \{0, \ldots, n_2-3\}$.
Suppose $n_1\leq n_2$, and put $m=n_2-n_1$
(recall that $N'(n_1, n_2, x, y)$ defined in the proof of Corollary \ref{t1p}).

Repeating the arguments form the proof of Corollary \ref{t1p}
we find out that
%If $x\geq y$, then
%$N(n_1, n_2, x, y)  = n_1+n_2+\min\{2y+n_1-x, x+y\} \leq
%n_1+n_2+\min\{2y+n_2-x+1, x+y\} = N(n_1, n_2, y, x)$. Thus further we
%assume that $x\leq y$.
%
%If $x<y$ and $x<n_1-3$, then
%$N(n_1, n_2, x, y)  = n_1+n_2+\min\{2x+n_2-y, x+y\} < N(n_1, n_2, x+1, y)$.
%Thus it remains to consider cases $x=y$ and $x=n_1-3$.
%
%We have $N(n_1, n_2, x, x)  = n_1+n_2+\min\{x+n_2, 2x\} \leq N(n_1, n_2, n_1-3, n_1-3)$.
%
it remains to maximize $N'(n_1, n_2, n_1-3, y)$ for $n_1-3\leq y\leq n_2-3$.

Let $y=n_1-3+z$, $z\in \{0, \ldots, m\}$.
We have $N'(n_1, n_2, n_1-3, y)  = n_1+n_2+\min\{2(n_1-3)+(n_1+m)-(n_1-3+z), 2(n_1-3)+z\} =
3n_1+n_2-6+\min\{m+3-z, z\}= 2n_1+2n_2-6 - l(m)$, where
$l(m) = m-\min\limits_{z\in \{0, \ldots, m\}} \{m+3-z, z\}$.
%Let $k(m) = \min\limits_{z\in \{0, \ldots, m\}} \{m-\min\{m+3-z, z\}\}$.
Note that $\min\{m+3-z, z\} \leq \left\lfloor \frac{m+3}{2} \right\rfloor$,
and  for $m\geq 2$ the equality holds for $z=\left\lfloor \frac{m+3}{2} \right\rfloor$
($\left\lfloor \frac{m+3}{2} \right\rfloor \leq m$ for $m\geq 2$);
thus for $m\geq 2$ we have $l(m) = m- \left\lfloor \frac{m+3}{2} \right\rfloor =
\left\lfloor \frac{m}{2} \right\rfloor -1.$
For $m=0$ and $m=1$ we have $l(0)=l(1)=0$.
From that the statement of Theorem \ref{t2p} follows.

%\newpage

%\section*{Список литературы}

%\input{references}

\end{document}